\documentclass[10pt,reqno]{amsart}

\usepackage{amsfonts}
\usepackage{amsmath}
\usepackage{amsthm}
\usepackage{amssymb}
\usepackage{latexsym}
\usepackage{multicol}
\usepackage{verbatim}
\usepackage{tabularx}
\usepackage{amscd,amssymb}

\advance\textwidth by 1.2in \advance\oddsidemargin by -.6in
\advance\evensidemargin by -.6in
\newtheorem*{cor}{Corollary}
\newtheorem*{lem}{Lemma}
\newtheorem*{prop}{Proposition}
\theoremstyle{definition}
\newtheorem*{defn}{Definition}
\theoremstyle{definition}
\newtheorem*{thm}{Theorem}

\newtheorem*{rem}{Remark}

\newenvironment{pf}{\proof}{\endproof}
\newcounter{cnt}
\newcounter{sbc}
\newcounter{tmp}
\newenvironment{enumerit}{\begin{list}{{\hfill\rm(\roman{cnt})\hfill}}{%
\settowidth{\labelwidth}{{\rm(viii)}}\leftmargin=\labelwidth%
\advance\leftmargin by
\labelsep\rightmargin=0pt\usecounter{cnt}}}{\end{list}}
\makeatletter
{\endlist}
\def\itlabel#1{\@bsphack\setcounter{tmp}{\@currentlabel}
  \protected@write\@auxout{}%
         {\string\newlabel{#1}{{\roman{tmp}}{\thepage}}}%
  \@esphack}
\makeatother
\theoremstyle{remark}

\numberwithin{equation}{section} 

\newtheorem{example}{Example}

\def\Hom{\operatorname{Hom}}
\def\End{\operatorname{End}}
\def\Ext{\operatorname{Ext}}

\def\Im{\operatorname{Im}}

\def\Ann{\operatorname{Ann}}

\def\ev{\operatorname{ev}}

\def\ord#1^#2{#1$^{\text{#2}}$}

\def\wdt#1{\widetilde{#1}}

\def\sop_#1^#2{\text{\scriptsize $\bigoplus\limits_{#1}^{#2}$}}

\newcommand{\thmref}[1]{Theorem~\ref{#1}}

\newcommand{\lemref}[1]{Lemma~\ref{#1}}
\newcommand{\propref}[1]{Proposition~\ref{#1}}
\newcommand{\corref}[1]{Corollary~\ref{#1}}

\newcommand{\id}{\operatorname{id}}
\newcommand{\ad}{\operatorname{ad}}
\newcommand{\tensor}{\otimes}
\newcommand{\nc}{\newcommand}
\newcommand{\rnc}{\renewcommand}

\nc{\cal}{\mathcal} \nc{\goth}{\mathfrak} \rnc{\bold}{\mathbf}

\renewcommand{\Bbb}{\mathbb}
\nc\bpi{{\mbox{\boldmath $\pi$}}} \nc\bvpi{{\mbox{\boldmath
$\varpi$}}}
 \nc\balpha{{\mbox{\boldmath $\alpha$}}}
\newcommand{\lie}[1]{\mathfrak{#1}}
\makeatletter
\def\section{\def\@secnumfont{\mdseries}\@startsection{section}{1}%
  \z@{.7\linespacing\@plus\linespacing}{.5\linespacing}%
  {\normalfont\scshape\centering}}
\def\subsection{\def\@secnumfont{\bfseries}\@startsection{subsection}{2}%
  {\parindent}{.5\linespacing\@plus.7\linespacing}{-.5em}%
  {\normalfont\bfseries}}
\makeatother
\newif\iffinal
\finaltrue
\def\subl#1{\subsection{\iffinal\else#1\fi}\label{#1}}

\nc{\Cal}{\cal} \nc{\Xp}[1]{X^+(#1)} \nc{\Xm}[1]{X^-(#1)}
\nc{\on}{\operatorname} \nc{\ch}{\mbox{ch}} \nc{\Z}{{\bold Z}}
\nc{\J}{{\cal J}} \nc{\C}{{\bold C}} \nc{\Q}{{\bold Q}}

\nc{\N}{{\Bbb N}} \nc\boa{\bold a} \nc\bob{\bold b}
\nc\boc{\bold
c} \nc\bod{\bold  d} \nc\boe{\bold e} \nc\bof{\bold f}
\nc\bog{\bold g} \nc\boh{\bold h} \nc\boi{\bold i} \nc\boj{\bold j} \nc\bok{\bold k} \nc\bol{\bold l} \nc\bom{\bold m}
\nc\bon{\bold n} \nc\boo{\bold o} \nc\bop{\bold p} \nc\boq{\bold
q} \nc\bor{\bold r} \nc\bos{\bold s} \nc\bou{\bold u}
\nc\bldx{\bold x}
\nc\boy{\bold y} \nc\bov{\bold v} \nc\bow{\bold w} \nc\boz{\bold
z}

\nc\ba{\bold A} \nc\bb{\bold B} \nc\bc{\bold C} \nc\bd{\bold D}
\nc\be{\bold E} \nc\bg{\bold G} \nc\bh{\bold H} \nc\bi{\bold I}
\nc\bj{\bold J} \nc\bk{\bold K} \nc\bl{\bold L} \nc\bm{\bold M}
\nc\bn{\bold N} \nc\bo{\bold O} \nc\bp{\bold P} \nc\bq{\bold Q}
\nc\br{\bold R} \nc\bs{\bold S} \nc\bt{\bold T} \nc\bu{\bold U}
\nc\bv{\bold V} \nc\bw{\bold W} \nc\bz{\bold Z} \nc\bx{\bold x}

\title[]{An application of free Lie algebras to current
algebras and their representation theory}
\author[]{Vyjayanthi Chari}
\address{Department of Mathematics, University of
California, Riverside, CA 92521, U.S.A.}
\email{chari@math.ucr.edu}
\author[]{Jacob Greenstein}\thanks{This work was partially supported by the Swiss National
Science Foundation}
\email{green@math.ucr.edu}

\begin{document}
\begin{abstract}
We realize the current algebra of a
Kac-Moody algebra as a quotient of a semi-direct product of the
Kac-Moody Lie algebra and the free Lie algebra of the Kac-Moody
algebra. We use this realization to study the representations of the current algebra.
In particular we see that every $\ad$-invariant ideal
in the symmetric algebra of the Kac-Moody algebra gives rise in a
canonical way to a representation of the current algebra. 
These  representations  include certain well-known families
of representations of the current algebra of a simple Lie algebra.
Another family of examples, which are the classical limits of the Kirillov-Reshetikhin modules, are also
obtained explicitly by using a construction of Kostant. Finally we study extensions in the category of
finite dimensional modules of the current algebra of a simple Lie algebra.
\end{abstract}
\maketitle

\section*{Introduction}
In this paper, we study the representation theory of current
algebras from a somewhat unusual viewpoint. Recall that the
current algebra $\lie a[t]$ of a Lie algebra $\lie a$ is the Lie
algebra $\lie a\otimes \bc[t]$ where the bracket is given in
the obvious  way, namely $[x\otimes f,y\otimes g]=[x,y]\otimes fg$
for all $x,y\in\lie a$, $f,g\in\bc[t]$. The algebra is naturally
graded by the non-negative integers  (so that the elements of
degree $k$ are of the form $\lie a\otimes t^k$) and has a unique
maximal graded ideal $\lie a[t]t:=\lie a\tensor t\bc[t]$.  The
starting point of our construction is to define a map from
the free Lie algebra $F(\lie a)$ on $\lie a$ to the maximal graded
ideal $\lie a[t]t$.  The free Lie algebra is a graded Lie algebra
where the grading is given by the positive integers and
whose enveloping algebra is the tensor algebra $T(\lie a)$ of
$\lie a$. Needless to say, the canonical Lie algebra homomorphism
$\tau$ from $F(\lie a)$ to $\lie a$ is not one of graded Lie
algebras. However, it is not hard to prove that if $\tau_k$ is the
restriction of $\tau$ to the $k^{th}$-graded piece, then the
direct sum $\lie k$ of the kernels of $\tau_k$ is a graded
ideal in $F(\lie a)$. Moreover, if we regard $T(\lie a)$
as an $\lie a$-module under the usual diagonal action induced by the
adjoint
action, we find that $F(\lie a)$ is an $\lie a$-submodule, $\lie k$ is a
$\lie a$-invariant ideal of
$T(\lie a)$ and that there exists a Lie algebra homomorphism
$\tau[t]: F(\lie a)\to\lie a[t]t$ with kernel $\lie k$. With a
little more work one then shows that this map actually extends to
a map from the natural semidirect product of $\lie a$ and $F(\lie a)$ 
to $\lie a[t]$. Moreover, the map is surjective if $\lie a$ is
its own commutator subalgebra, which we shall assume from now on.

Using this realization, we then see in particular that to any
representation $V$ of $\lie a[t]$ one can associate a $\lie a$-invariant 
ideal $\bi(V)$ of $T(\lie a)$ containing $\lie k$.
Conversely, any  $\lie a$-invariant  ideal $\bi$ of $T(\lie a)$
containing $\lie k$ defines a representation $V(\bi)$ of $\lie a[t]$ 
such that $\bi(V(\bi))=\bi$. However, non-isomorphic modules
can give rise to the same ideal. Thus for example, any module $V$
which is obtained from a module for $\lie a$ by evaluating at zero
corresponds to the augmentation ideal in $T(\lie a)$. Nonetheless,
since there are infinitely many $\lie a$-invariant ideals
$\bi\supset\lie k$ of $T(\lie a)$ we get a large family of
explicitly defined representations of $\lie a[t]$. For instance,
by observing that $\lie k$ is contained in the two sided ideal of
$T(\lie a)$ generated by elements of the from $x\otimes y-y\otimes x$, 
$x,y\in\lie a$, we find that  every  $\lie a$-invariant ideal
in $S(\lie a)$  corresponds to a $\lie a[t]$-module $M$ which
satisfies $(x\otimes t^2)M=0$ for all $x\in\lie a$. The study of such
ideals in $S(\lie a)$ is well-established, at least for $\lie a$ 
semi-simple (see for example~\cite{Bo,Jo,McG}) and we expect that
the corresponding modules for $\lie a[t]$ should be similarly
interesting. Indeed, one can see that a subfamily of the classical
counterparts of the Kirillov-Reshetikhin modules can be described
in this language. The results of Section 2 of this paper are
somewhat more general than the above outline and we use a
construction of Kostant that also allows us to construct the
fundamental Weyl modules for the current algebra (cf.~\cite{Ckirres}, \cite{CPweyl}).

In Section~3, we use the construction to describe extensions
between $\lie a[t]$-modules which split as $\lie a$-modules. We
then use this to describe extensions between irreducible
finite-dimensional modules for the current algebra 
associated to a simple finite dimensional Lie algebra.

\subsection*{Acknowledgments}We
are grateful to A.~Alekseev, A.~Berenstein, A.~Joseph and E.~Mukhin for
stimulating discussions.

\section{A realization of the current algebra}
Let $\lie a$ be a complex Lie algebra satisfying  $[\lie a,\lie
a]=\lie a$. In this section we first realize the maximal ideal in
the current algebra $\lie a[t]$  as a quotient of the free Lie
algebra $F(\lie a)$ on the underlying vector space $\lie a$.
Further, since $\lie a$ acts in a natural way on $F(\lie a)$ we
see that the current algebra is a quotient of  $\lie F$ which is
the semi-direct product of $\lie a$ and $F(\lie a)$. Although we
work over the complex numbers, since that is our primary interest
in this paper, these results remain valid over any ground field.

\subl{Map10}  Let $\bz_+$ (respectively, $\bn$) be the set of non-negative
(resp. positive) integers. Let $\lie a$ be a complex
 Lie algebra with $\lie a=[\lie a,\lie a]$,  $U(\lie a)$ the universal
enveloping  algebra of $\lie a$ and $T(\lie a)$ the tensor algebra
of $\lie a$. Clearly, $T(\lie a)=\bigoplus_{r\ge 0} T^r(\lie a)$
is a $\bz_+$-graded associative algebra.  Set
$$
F^1(\lie a)=T^1(\lie a)=\lie a
$$
and let  $F^r(\lie a)$ be  the subspace of $T^r(\lie a)$ spanned
by elements of the form $x\tensor\boy-\boy\tensor x$, $x\in\lie
a$, $\boy\in F^{r-1}(\lie a)$. Then
$$
F(\lie a)=\bigoplus_{r> 0} F^r(\lie a),
$$
is the free Lie algebra on $\lie a$ with the bracket operation
given by
$$[\boy,\boy']_T=\boy\otimes\boy'-\boy'\otimes\boy,\qquad
\boy,\boy'\in F(\lie a). $$ Obviously $F(\lie a)$ is a
$\bn$-graded Lie algebra and  $T(\lie a)$ is the universal
enveloping algebra of $F(\lie a)$. Let $\tau: F(\lie a)\to \lie a$
be the canonical homomorphism of Lie algebras and we also denote
by $\tau$ the extension to the map of  associative algebras
$T(\lie a)\to U(\lie a)$.
 Let  $\tau_r:F^r(\lie a)\to\lie a$
be the restriction of $\tau$ to $F^r(\lie a)$. Note that
$\tau_1=\id$ and for all $r,s\in\bn_+$,
\begin{equation}\label{lietau}
[\tau_r(\bx),\tau_s(\boy)]_{\lie a}=\tau_{r+s}([\bx,\boy]_T),\qquad \bx\in
F^r(\lie a), \boy\in F^s(\lie a).
\end{equation} Let $\ad_T:\lie a\to \End(T(\lie a))$ be the usual diagonal action of $\lie a$ on
$T(\lie a)$ induced by $\ad:\lie a\to \End(\lie a)$, namely
$$\ad_T(x)(y_1\otimes\cdots\otimes
y_r)=\sum_{i=1}^ry_1\otimes\cdots\otimes
[x,y_i]\otimes\cdots\otimes y_r,$$ for all $x,y_1,\cdots
,y_r\in\lie a$, $r\in\bn$.

\subl{Map15} Let $\bc[t]$ be the ring of polynomials in an
indeterminate $t$. The space $\lie a[t]=\lie a\otimes \bc[t]$ is
naturally a Lie algebra with the bracket defined by 
$$ 
[x\otimes f,y\otimes g]_{\lie a[t]}=[x,y]_{\lie a}\otimes fg,\quad x,y\in\lie a, \quad
f,g\in\bc[t], $$ and let $\lie a[t]t$ denote the ideal $\lie
a\otimes t\bc[t]$ of $\lie a[t]$. Henceforth we write $x f$ for
the element $x\tensor f$, $x\in\lie a$, $f\in\bc[t]$ of $\lie
a[t]$.
\begin{thm}
\begin{enumerit}
\item\itlabel{Map15.i} The map $\tau[t]: F(\lie a)\to \lie a[t]$
defined by $ \tau[t](\bldx)=\tau(\bldx)t^r$, for $\bldx\in
F^r(\lie a)$ and $ r\ge 1 $ is a homomorphism of graded Lie
algebras and $$ \ker\tau[t]=\bigoplus_{r> 0}\ker\tau_r. $$
Moreover, $\ker\tau[t]$ is contained in the two-sided ideal of
$T(\lie a)$ generated by elements $x\otimes y-{y\otimes x}-[x,y]$ for all
$x,y\in\lie a$.
\item\itlabel{Map15.ii} Let $\lie{F}=\lie a\oplus
F(\lie a)$ and set
$$
[x,y]_{\lie F}=[x,y]_{\lie a}, \quad [x,\boy]_{\lie
F}=\ad_T(x)\boy=-[\boy,x]_{\lie F}, \quad [\boy,\boz]_{\lie
F}=[\boy,\boz]_T
$$
for all~$x,y\in\lie a$, $\boy,\boz\in F(\lie a)$. Then~$\lie F$ is
a Lie algebra and
$$
\lie a[t]\cong  \lie F/\ker\tau[t],
$$
where $\tau[t]$ is extended to~$\lie F$ by
setting~$\tau[t](x)=x\tensor 1$, $x\in\lie a$.
\end{enumerit}
\end{thm}
\begin{pf} 
Take $\bx\in F^r(\lie a)$, $\boy\in F^s(\lie a)$. Then $[\bx,\boy]_T\in F^{r+s}(\lie a)$
and so by~\eqref{lietau}
$$
\tau[t]([\bx,\boy]_T)=\tau_{r+s}([\bx,\boy]_T)t^{r+s}=([\tau_r(\bx),\tau_s(\boy)]_{\lie a})t^{r+s}=
[\tau_r(\bx)t^s,\tau_s(\boy)t^s]_{\lie a[t]}=[\tau[t](\bx),\tau[t](\boy)]_{\lie a[t]},
$$
which proves the first statement in~\eqref{Map15.i}.
The second statement is
obvious. To prove that $\lie F$ is a Lie algebra,  note that for
$x,y\in\lie a$, $\boz\in F(\lie a)$,  we have
\begin{equation}\label{15}
\ad_T(x)[y,\boz]_T=[[x,y],\boz]_T+[y,\ad_T(x)(\boz)]_T.
\end{equation}
Since~$F^r(\lie
a)$ is spanned by elements~$[y,\boz]_T$, $y\in\lie a$, $\boz\in
F^{r-1}(\lie a)$, an obvious induction on~$r$ proves that
$\ad_T(x)(F^r(\lie a))\subset F^r(\lie a)$ for all~$r\in\bn$. This implies that
$$
\ad_T(x)[\boy,\boz]_T=[\ad_T(x)\boy,\boz]_T+[\boy,\ad_T(x)\boz]_T
$$
for all $x\in\lie a$, $\boy,\boz\in F(\lie a)$, whence
$$
[x,[\boy,\boz]_{\lie F}]_{\lie F}+[\boz,[x,\boy]_{\lie F}]_{\lie F}+
[\boy,[\boz,x]_{\lie F}]_{\lie F}=0.
$$
Since $F^r(\lie g)$ is preserved by~$\ad_T$ by the above, the Jacobi identity for
$x,y\in\lie a$, $\boy\in F(\lie a)$ is trivially checked. Finally, to prove that the extension of
$\tau[t]$ to $\lie F$ is a homomorphism of Lie algebras, it is
enough to check that for all $r\in\bn$, the map $\tau_r:F^r(\lie
a)\to \lie a$ is a surjective map of $\lie a$-modules. This is
also proved by an induction on $r$. Since $\tau_1=\id$, the induction starts. Applying~$\tau$ to both sides
of~\eqref{Map15} with~$\boz\in F^{r-1}(\lie g)$, we get
\begin{align*}
\tau(\ad_T(x)([y,\boz]_T))&=
[[x,y],\tau(\boz)]-[\tau(\ad_T(x)\boz),y]=
[[x,y],\tau(\boz)]-[[x,\tau(\boz)],y]\\&=[x,[y,\tau(\boz)]]
=[x,\tau([y,\boz]_T)]
\end{align*}
since~$\tau:F^{r-1}(\lie a)\to\lie a$ is a~$\lie a$-module
homomorphism by the induction hypothesis. It remains to prove
that~$\tau:F^r(\lie a)\to\lie a$ is surjective. Since $\lie
a=[\lie a,\lie a]$, for any~$x\in\lie a$ there exist $y,z\in\lie
a$ such that~$x=[y,z]$. By the induction hypothesis, there exists
$\boz\in F^{r-1}(\lie a)$ such that~$z=\tau(\boz)$.
Then~$x=[y,\tau(\boz)]=\tau([y,\boz]_T)$,
whence~$x\in\tau(F^r(\lie a))=\Im\tau_r$.
\end{pf}
\begin{cor} For all $k\ge 0$, we have an isomorphism of Lie algebras $$F(\lie
a)\left/\Big(\ker\tau[t]+\sum_{r>k}F^r(\lie a)\Big)\right.\cong
\lie a[t]/\lie a\otimes t^k\bc[t].
$$
\end{cor}

\section{Construction of $\lie a[t]$-modules}

In this section we study some applications of the realization of
$\lie a[t]$ provided by~\thmref{Map15} to the representation theory of $\lie a[t]$.
 Thus, if  $\rho[t]:\lie a[t]\to \End (V)$ is a homomorphism
of Lie algebras, it is clear that the action of $\lie a[t]$ on $V$
is determined by  the images of $\rho[t](x)$ and~$\rho[t](xt)$
in~$\End(V)$ for all~$x\in\lie a$. In this section we give a
necessary and sufficient condition for a pair of linear maps from
$\lie a$ to $\End(V)$ to define a $\lie a[t]$-module structure
on~$V$ and describe the isomorphism classes of such structures.

\subl{Cons1} Given a collection of vector spaces~$V_j$,
$j=1,\dots,k+1$ and maps~$\xi_j\in\Hom(\lie a,\Hom(V_j,V_{j+1}))$,
$j=1,\dots,k$ define $\xi_k\odot\cdots\odot\xi_1\in\Hom(T^k(\lie
a),\Hom(V_1,V_{k+1}))$ by
$$ (\xi_k\odot\cdots
\odot\xi_1)(x_k\tensor\cdots \tensor x_1):= \xi_k(x_k)\circ\cdots
\circ \xi_1(x_1).
$$
Let $\xi:\lie a\to \cal A$ be a map of vector spaces from $\lie a$
to an associative algebra with unity $\cal A$. Denote by
$\xi_T:T(\lie a)\to\cal A $ the natural algebra homomorphism
extending $\xi$. In particular, $\xi_T|_{T^0(\lie a)}$ is defined
by $\xi_T(c)=c1_{\cal A}$, $c\in T^0(\lie a)=\bc$ and
$$
\xi_T|_{T^r(\lie a)}=\xi^{\odot r}.
$$
Let  $\xi_F$ be
the restriction of $\xi_T$ to $ F(\lie a)$. Clearly, $\xi_F$ is a
homomorphism of Lie algebras, where $\cal A$ is regarded as a Lie algebra 
in a natural way.
\begin{lem}
If  $\xi:\lie a\to \cal A$ is a homomorphism of Lie algebras
then
\begin{equation}\label{rdot} \xi^{\odot r}|_{F^r(\lie a)}=\xi\circ
\tau_r,\qquad r>0.
\end{equation}
\end{lem}
\begin{pf}
The proof is by induction on~$r$. If $r=1$, then $\tau_1$ is the
identity map and hence induction begins. For the inductive step,
it is sufficient to prove the assertion for $\bldx=[x,\boy]_T$
where $x\in\lie a$, $\boy\in F^r(\lie a)$. Then
\begin{equation*}
\xi^{\odot r+1}(\bldx)= [\xi(x),\xi^{\odot
r}(\boy)]_{\cal A}=[\xi(x),\xi(\tau_r(\boy))]_{\cal A}=
\xi([x,\tau_r(\boy)])=\xi(\tau_{r+1}(\bldx)).\qedhere
\end{equation*}
\end{pf}

\subl{Cons10} Suppose that we are given a homomorphism of Lie
algebras $\rho:\lie a\to \End(V)$, that is a $\lie a$-module
structure on $V$. Then $\End(V)$ is an $\lie a$-module in a
natural way. Suppose we are given $\eta:\lie a\to \End(V)$ a map
of $\lie a$-modules, that is
$$
[\rho(x),\eta(y)]=\eta([x,y]),\qquad x,y\in\lie a.
$$
Then $\eta_F$ is a map of $\lie a$-modules and we have
$$
[\rho(x),\eta^{\odot r}(\boy)]=\eta^{\odot r}(\ad_T(x)\boy),\qquad
x\in\lie a,\,\boy\in F^r(\lie a).
$$
\begin{prop} Let $\rho:\lie a\to\End(V)$ be a
homomorphism of Lie algebras. There exists a homomorphism 
$\rho[t]:\lie a[t]\to\End(V)$ of Lie algebras such that
$\rho[t](x\otimes 1)=\rho(x)$ for all $x\in\lie a$, if and only if
there exists $\eta\in\Hom_{\lie a}(\lie a,\End(V))$ such that
$$\ker\tau[t]\subset\ker \eta_T ,$$  where we assume that the
$\lie a$-module structure on  $\End(V)$ is the natural one given
by $\rho$.
\end{prop}
\begin{pf}
 Suppose first that $\rho[t]$ exists. Define $\eta:\lie a\to\End(V)$ by $\eta(x)=\rho[t](xt)$. Since $$
\eta([x,y])(v)=([x,y]t)v=([x,yt])v=\rho(x)(\eta(y)(v))-\eta(y)(\rho(x)(v)),
$$ it follows that $\eta$ is a map of $\lie a$-modules. To prove
that $\eta_T(\ker\tau_r)=0$ it suffices to prove that
\begin{equation}\label{Cons10.10}
\eta^{\odot r}(\bldx)(v)=(\tau_r(\bldx)t^r) v\qquad\text{for
all~$\bldx\in F^r( \lie a)$ and for all~$r>0$}.
\end{equation}
We proceed by induction on $r$, the case~$r=1$ being obvious.
For the inductive step, it is sufficient to
prove~\eqref{Cons10.10} for all elements of~$F^{r+1}(\lie a)$ of
the form $[x,\boy]_T$, $x\in\lie a$, $\boy\in F^r(\lie a)$. One
has
\begin{align*}
\eta^{\odot r+1}([x,\boy]_T)v&=[\eta(x),\eta^{\odot
r} (\boy)]v\\ &=(xt)(\tau_r(\boy)t^r)
v- (\tau_r(\boy)t^r)(xt)v\\
&=([xt,\tau_r(\boy)t^r])v=([x,\tau_r(\boy)]t^{r+1})v
\\
&=(\tau_{r+1}([x,\boy]_T)t^{r+1})v,
\end{align*}
where the second equality is a consequence of the induction
hypothesis.

For the converse, suppose that we are given a map $\eta:\lie
a\to\End(V)$ of $\lie a$-modules and let $\eta_F: F(\lie a)\to
\End(V)$ be the corresponding  homomorphism of Lie algebras, with
$\eta_F(\bx)=\eta^{\odot r}(\bx)$ if $\bx\in F^r(\lie a)$. Since
$\eta_T(\ker\tau_r)=0$ for all $r\ge 1$, it follows immediately
that we have a homomorphism of Lie algebras, $\eta[t]: \lie
a[t]t\to \End(V)$. Moreover, since $\eta$ is a homomorphism of
$\lie a$-modules, the map $\eta_F$ and hence the induced map
$\eta[t]$ are $\lie a$-module homomorphisms. It is now sufficient
to prove that the formulas
$$
\rho[t](x)=\rho(x),\quad \rho[t](x)=\eta[t](xt^r),\qquad x\in\lie
a,\, r\ge 1,
$$
define a representation of $\lie a[t]$. For this, it suffices to
show that
$$[\rho[t](x),\rho[t](yt^r)]=\rho[t]([x,yt^r]).
$$
But this equivalent
to proving that
$$[\rho(x),\eta[t](yt^r)]=\eta[t]([x,yt^r]) $$
which is the statement that $\eta[t]$ is a map of $\lie
a$-modules.
\end{pf}

\subl{Cons20} Motivated
by~\propref{Cons10} we introduce the following definition.
\begin{defn}
Let $\cal H(V)$ be the subset of $\Hom(\lie a,\End(V))\times
\Hom(\lie a,\End(V))$ consisting of pairs $(\rho,\eta)$ which
satisfy,
\begin{enumerit}
\item[$(C_1)$] $\rho$ is a homomorphism of Lie algebras.
\item[$(C_2)$] $\eta$ is a homomorphism of $\lie a$-modules where
$\End(V)$ acquires a $\lie a$-module structure through~$\rho$.
\item[$(C_3)$] For all~$r>0$, $\ker\tau_r\subset\ker\eta_T$.
\end{enumerit}
\end{defn}
\noindent Given a homomorphism of Lie algebras $\rho:\lie
a\to\End(V)$, set $$ \cal H_\rho(V)=\{\eta\in\Hom_{\lie a}(\lie
a,\End(V)): (\rho,\eta)\in\cal H(V)\}. $$ Given a pair
$(\rho,\eta)\in \cal H(V)$, denote by $V(\rho,\eta)$ the $\lie
a[t]$-module structure defined on~$V$ by~$\rho$, $\eta$ as
in~\propref{Cons10}. Note that $\eta=0$ is in $\cal H_\rho(V)$.
For $a\in\bc$, let $\epsilon_a$ be the endomorphism of $\lie a[t]$
defined by $\epsilon_a(xt^r)=a^r xt^r$ and $\gamma_a$ be the
automorphism of $\lie a[t]$ defined by $\gamma_a(xt^r)=x(t-a)^r$.
\begin{lem} Let $\eta\in \cal H_\rho(V)$. Then,
 for all $a\in\bc$, we have
\begin{enumerit}
\item\itlabel{Cons20.i} $a\eta\in \cal H_\rho(V)$,
\item\itlabel{Cons20.ii} $\eta-a\rho\in\cal{H}_\rho (V)$,
\item\itlabel{Cons20.iii}  $a\rho\in\cal H_\rho(V)$.
\end{enumerit}
\end{lem}
\begin{pf}  To prove \eqref{Cons20.i} (respectively, \eqref{Cons20.ii}) it suffices to observe that
$V(\rho,a\eta)$ (resp. $V(\rho,\eta-a\rho)$) is isomorphic to the
pull back through $\epsilon_a$ (resp. $\gamma_a$) of $V(\rho,\eta)$. Part
\eqref{Cons20.iii} follows immediately from~\eqref{Cons20.ii} by taking $\eta=0$.
\end{pf}
\begin{rem} A $\lie a[t]$-module of the form $V(\rho,a\rho)$, $a\in\bc$
is an evaluation module, namely it is obtained by pulling the
$\lie a$-module $V$ back through the evaluation homomorphism $\ev_a:\lie a[t]\to\lie a$
defined by $\ev_a(x f)=f(a)x$ for all $x\in\lie g$, $f\in\bc[t]$. It is tempting to
conjecture that if $\rho$ is irreducible then $\cal H_\rho(V)=\bc\rho$.
In fact, when $\lie g$ is a simple Lie algebra and $V$ is an
irreducible highest weight representation of $\lie g$, one can
prove without too much difficulty  that this statement is true.
However, the general case is far from clear.
\end{rem}

\subl{Cons22} The following is trivial.
\begin{lem}
A map $\phi\in\Hom_{\lie a}(V(\rho,\eta),V'(\rho',\eta'))$ is a
homomorphism of $\lie a[t]$-modules if and only if
$\phi\circ\eta(x)=\eta'(x)\circ\phi$ for all $x\in\lie a$. In
particular, given $\eta,\eta'\in\cal H_\rho(V)$ the $\lie
a[t]$-modules  $V(\rho,\eta)$ and $V(\rho,\eta')$ are isomorphic
if and only if there exists a $\lie a$-module automorphism $\phi$
of $V$ such that $\eta'(x)=\phi\circ\eta(x)\circ\phi^{-1}$ for all
$x\in\lie a$.
\end{lem}

\subl{Cons25}

In the rest of this section, we  use the above construction to
produce a large family of examples of non-isomorphic $\lie
a[t]$-modules $V(\rho,\eta)$.

For $\eta\in\cal H_\rho(V)$ set
$$
\bi(\rho,\eta)=\ker\eta_T.
$$
Clearly, $\bi(\rho,\eta)$ is a two-sided  $\ad_T$-invariant ideal
of $T(\lie a)$ containing $\ker\tau[t]$.  It follows immediately
from~\lemref{Cons22} that $\bi(\rho_1,\eta_1)=\bi(\rho_2,\eta_2)$
if $V_1(\rho_1,\eta_1)\cong V_2(\rho_2,\eta_2)$. Observe that for
all $a\in\bc^\times $ we have $\bi(\rho, a\eta)=\bi(\rho,\eta)$
for all $a\not=0$. Moreover, $\bi(\rho,0)$ is the augmentation
ideal $T(\lie a)_+$ of $T(\lie a)$. In the case when $\eta=a\rho$
for $a\in\bc^\times$, we have that   $\eta_T:T(\lie a)\to \End(V)$
factors through to an associative algebra homomorphism
$\eta_U:U(\lie a)\to\End(V)$. Obviously,
$\ker\eta_U=\operatorname{Ann}_{U(\lie a)} V$. It follows that
$\bi(\rho,a\rho)$ is  the inverse image of $\Ann_{U(\lie a)} V$
with respect to $\tau:T(\lie a)\to U(\lie a)$.

Conversely, given any two sided $\lie a$-invariant ideal in
$T(\lie a)$ which contains $\ker\tau[t]$, we can define a
representation $V(\rho_\bi,\eta_\bi)$ as follows. Set $V=T(\lie
a)/\bi$. Then the action $\ad_T$ on $T(\lie a)$ induces a  natural
$\lie a$-module structure on $V$ and we denote $\rho_\bi$ the
corresponding homomorphism of Lie algebras $\lie a\to \End(V)$.
Define $\eta_\bi: \lie a\to\End(V)$ by
$$
\eta_\bi(x)(\boy+\bi)=x\otimes\boy+\bi, \quad x\in\lie a,\,
\boy\in T(\lie a)
$$

Let $\lie M$ be the set of of  two sided $\ad_T$-invariant ideals
in $T(\lie a)$ containing $\ker\tau[t]$. Denote by $\cal M$ the
set of isomorphism classes of $\lie a[t]$-modules.
\begin{prop}
\begin{enumerit}
\item\itlabel{Cons25.0} $\eta_\bi\in\cal H_{\rho_\bi}(T(\lie
a)/\bi)$.
\item\itlabel{Cons25.i} $\bi(\rho_\bi,\eta_\bi)=\bi$.
\item\itlabel{Cons25.ii} The map $\lie M\to \cal M$ defined by
$\bi\mapsto [V(\rho_\bi,\eta_\bi)]$ is injective
\item\itlabel{Cons25.iii} The map $\cal M\to \lie M$ defined by
$[V(\rho,\eta)]\mapsto \bi(\rho,\eta)$ is surjective.
\end{enumerit}
\end{prop}
\begin{pf}
For~\eqref{Cons25.0}, we should prove first that $\eta_\bi$ is a
homomorphism of $\lie a$-modules. Indeed, for all $x,y\in\lie a$,
$\boz\in T(\lie a)$ we have $$
(\rho_\bi(x)\eta_\bi(y)-\eta_\bi(y)\rho_\bi(x))(\boz+\bi)=\rho_\bi(x)(y\tensor
\boz)-y\tensor \rho_\bi(x)\boz+\bi =[x,y]\tensor
\boz+\bi=\eta_\bi([x,y])(\boz+\bi). $$ Furthermore,
$\bi\subset\ker(\eta_\bi)_T$ by definition and so~$\eta_\bi\in\cal
H_{\rho_\bi}(T(\lie a)/\bi)$. In order to prove~\eqref{Cons25.i}
it remains to show that $\ker (\eta_\bi)_T\subset\bi$. Since
$\bx\in\ker(\eta_\bi)_T$ implies that $\bx\tensor\boy\in\bi$ for
all $\boy\in T(\lie a)$, the result  follows by taking $\boy=1$.
Finally, \eqref{Cons25.ii} and~\eqref{Cons25.iii} follow
immediately from~\eqref{Cons25.i}.
\end{pf}
 It is easy to
see that the ideal  $T(\lie a)_+$ gives rise to  the trivial
module for $\lie a[t]$. Furthermore, consider the ideal
$\bi=\ker\tau$ in $T(\lie a)$. Since $T(\lie a)/\ker\tau\cong
U(\lie a)$, it follows that $\ker\tau$ is the two-sided ideal in
$T(\lie a)$ generated by $x\tensor y-y\tensor x-[x,y]$ for all
$x,y\in\lie a$. Then it is easy to check that $\ker\tau$ is an
$\ad_T$-invariant ideal and
$\ker\tau[t]\subset\ker\tau$ by~\thmref{Map15}\eqref{Map15.i}. The $\lie a[t]$-module corresponding
to this ideal by the above construction identifies with
$U(\lie a)(\rho,\eta)$ with $\rho(x)u=\ad(x)u$ and $\eta(x)u=xu$.
It follows that $\eta^{\odot r}(\bx)u=\tau(\bx)u$ for all $\bx\in
F^r(\lie a)$. In particular, in this module one has $(xt^r)u=xu$
for all $r>0$.

\subl{Cons35} Consider now the special case of $\lie a[t]$-modules
$V(\rho,\eta)$ satisfying the condition $xt^r V(\rho,\eta)=0$ for
all $x\in\lie a$, $r\ge 2$.  Denote by~$\cal M_2$ the set of
isomorphism classes of such $\lie a[t]$-modules. By~\corref{Map10}
the isomorphism class of $V(\rho,\eta)$ is in $\cal M_2$ if and
only if $\eta_F(F^2(\lie a))=0$, that is $[\eta(x),\eta(y)]=0$ for
all $x,y\in\lie a$. It follows that $\eta_T$ factors through to an
algebra homomorphism $\eta_S: S(\lie a)\to\End(V(\rho,\eta))$
which is also a $\lie a$-module map. Let
$I(\rho,\eta)=\ker\eta_S$. Then $I(\rho,\eta)$ is a $\lie
a$-invariant ideal of $S(\lie a)$. Denote the set of such ideals
by $\lie M_S$.

On the other hand, given $I\in\lie M_S$, observe that $V=S(\lie
a)/I$ is a $\lie a$-module in a natural way. Denote the
corresponding Lie algebra homomorphism $\lie a\to\End(V)$ by
$\rho_I$ and define $\eta_I$ by
$$
\eta_I(x)(y+I)=xy+I,\qquad \forall\,x\in\lie a,\,y\in S(\lie a).
$$
\begin{prop}
\begin{enumerit}
\item\itlabel{Cons35.0} $\eta_I\in\cal H_{\rho_I}(S(\lie a)/I)$.
\item\itlabel{Cons35.i} $I(\rho_I,\eta_I)=I$.
\item\itlabel{Cons35.ii} The map $\lie M_S\to \cal M_2$ defined by
$I\mapsto [V(\rho_I,\eta_I)]$ is injective
\item\itlabel{Cons35.iii} The map $\cal M_2\to \lie M_S$ defined
by $[V(\rho,\eta)]\mapsto I(\rho,\eta)$ is surjective.
\end{enumerit}
\end{prop}
The proof of this proposition is similar to that
of~\propref{Cons25} and is omitted. Notice that if $I=0$, then
$S(\lie a)$ is a $\lie a[t]$-module with $\eta(x)$ given by
$$
\eta(x)y=xy.
$$
In particular, $S_k=\sum_{j\ge k} S^j(\lie a)$ is a $\lie
a[t]$-submodule of $S(\lie a)$ and the quotient $\lie
a[t]$-module $V_k=S_{k-1}/S_k$ is an evaluation module at $a=0$.
If $\lie a$ admits a non-degenerate symmetric invariant bilinear form
$(\cdot,\cdot)_{\lie a}$ then,
passing to the dual module, one obtains another $\lie a[t]$-module
structure on $S(\lie a)$ with $\eta$ given by
$$
\eta(x)(x_1\cdots x_k)=\sum_{j=1}^k (x,x_j)_{\lie a} x_1\cdots \hat
x_j\cdots x_k.
$$
The corresponding ideal is again~$I=0$. Next, consider the ideal
$I\subset S(\lie a)$ generated by $S^2(\lie a)$. Then as $\lie
a$-modules, we have
$$
V(\rho_I,\eta_I)=\bc\oplus \lie a,
$$ and
$\eta_I(x)(a,y)=(0,ax)$ for all $a\in\bc$, $x,y\in\lie a$. In the
case of simple Lie algebras this example can be further
generalized and which we now discuss.

\subl{Cons50}  Assume that $\lie g$ is a simple finite-dimensional
complex simple Lie algebra. Let $\theta$ be the highest root of
$\lie g$ and let $x^+_\theta$ be a non-zero element in the
corresponding root space. 
Set $x^+_\theta(r)=(x^+_\theta)^{\otimes r}\in T^r(\lie g)$. 
It is well-known that the $\lie g$-module $V(r\theta)= U(\lie
g)x^+_\theta(r)$ is an irreducible submodule of $T^r(\lie g)$ and
moreover that we can write
$$
T^r(\lie g)=V(r\theta)\oplus \tilde K_r,
$$
where $\tilde K_r$ is a $\lie g$-module. Further, it is not hard to see that if $x\in\lie g$
and $\boy\in \tilde K_r$, then $x\otimes \boy, \boy\otimes x\in\tilde K_{r+1}$ and so $\tilde K=\bigoplus_{r\ge 1} \tilde
K_r$ is a two-sided $\ad_T$-invariant ideal in $T(\lie g)$.
We claim that the quotient
$$
T(\lie g)/\tilde K=\bigoplus_{r\ge 0} V(r\theta)
$$
is a $\lie g[t]$-module and, moreover, its isomorphism class is in~$\cal M_2$.
For, notice that
$$
V(2\theta)\subset\langle\{x\otimes y+y\otimes x: x,y\in\lie g\}\rangle,
$$
and hence the elements of the form $x\otimes y-y\otimes x\in
\tilde K_2$ for all $x,y\in\lie g$. This implies that the ideal
generated by such elements is in the ideal $\tilde K$. In particular,
$F^r(\lie g)\subset\tilde K$ for all $r>1$, which
proves our claim.
Moreover, since
$$
\eta_{\tilde
K}(x^+_\theta)(x^+_\theta(r))=x^+_\theta(r+1),
$$
it follows easily that this module is indecomposable.

For $r\ge 1$, let  $K_r$  be the ideal
generated by $\tilde K$ and $T^r(\lie g)$.
This gives rise to a family of indecomposable $\lie g[t]$-modules
$$
T(\lie g)/K_r\cong\bigoplus_{s=0}^{r-1} V(s\theta).
$$
These modules (or rather their duals) appear in the literature (cf.~\cite{Ckirres},\cite{KR}) 
and are known to be the classical limits
of a family of  Kirillov-Reshetikhin modules.

\subl{Cons60} We now use a construction due to Kostant (cf.~\cite{K}) to give further
examples of modules in $\cal M_2$ for a simple Lie algebra $\lie
g$.

 To describe this construction, assume that   $\rho:\lie g\to
 \End(V)$ is a self-dual finite-dimensional representation of $\lie g$.
 Let $(\cdot,\cdot)$ be a non-degenerate symmetric $\lie
g$-invariant bilinear form on $V$. Let $\lie{so}(V)$ be the
corresponding Lie algebra of skew-symmetric endomorphisms of $V$
and notice that $\Im\rho\ \subset\lie{so}(V)$. Define
$\zeta:\bigwedge^2 V\to\lie{so}(V)$ by extending linearly
$$
\zeta(u_1\wedge u_2)v=(u_1,v)u_2-(u_2,v)u_1,  \qquad
u_1,u_2, v\in V.
$$
It is easy to see that $\zeta$ is injective and hence an
isomorphism of vector spaces. Let $\varphi:\lie g\to \bigwedge^2
V$ be the map ~$$\varphi(x)=\zeta^{-1}\circ \rho(x), \ x\in\lie
g.$$ The exterior powers $\bigwedge ^kV$, $k\ge 0$   are all $\lie
g$-modules under the usual diagonal action and one can show that
$\varphi$ is a homomorphism of $\lie g$-modules. Moreover, if we
set
$$
\eta(x)(u)=\varphi(x)\wedge u,\qquad x\in\lie g,\,
u\in\textstyle\bigwedge V,
$$ it is easy to see that the pair
$(\wedge\rho,\eta)$ defines a $\lie g[t]$-module structure on
$\End(\bigwedge V)$, here $\wedge \rho$ is the diagonal action of
$\lie g$ on $\bigwedge V=\bigoplus_{k\ge 0}\bigwedge^k V$.
Since $u\wedge u'=u'\wedge u$ for all $u,u'\in\bigwedge^2 V$,
it follows that $[\eta(x),\eta(y)]=0$ for all $x,y\in\lie g$.

For the dual module $ (\wedge^*\rho,\eta^*)$ we obtain the
following formula for~$\eta^*(x)$,
$$
\eta^*(x)(v_{1}\wedge\cdots\wedge v_{k})= -\sum_{1\le r<s\le k}
(-1)^{r+s} (\varphi(x),v_{r}\wedge v_{s})_{V}v_{1}\wedge \cdots \wedge
\hat v_{r}\wedge\cdots \wedge \hat v_{s}\wedge \cdots\wedge v_{k},
$$ where~$v_{r}\in V$ for $1\le r\le k$ and the form $B_V$ is extended
to $\bigwedge V$ by setting $B_V(v_1\wedge\cdots \wedge v_k,
w_1\wedge\cdots\wedge v_r)=\delta_{r,k}\det((B_V(v_i,w_j))_{i,j=1}^k)$.

In the special case, when $\lie g=\lie{so}(n)$ and $V$ is the
natural representation of $\lie g$,  it is not hard to see that
the subspaces $\bigoplus_{k\ge r}\bigwedge^{2k} V$ and $\bigoplus_{k\ge
r}\bigwedge^{2k+1} V$  are indecomposable $\lie g[t]$-submodules of $\bigwedge V(\wedge \rho,\eta)$.
Similarly, the dual modules $ (\wedge^*\rho,\eta^*)$ give rise to
a family of indecomposable modules $\bigoplus_{0\le k\le
r}\bigwedge^{2k} V$ and $\bigoplus_{0\le k\le r}\bigwedge^{2k+1} V$.
These modules are the classical limits of the fundamental
Kirillov-Reshetikhin modules or Weyl modules (cf.~\cite{CPweyl}).

\section{Extensions}

The main result of this section is a description of  the vector
space of extensions between two $\lie a[t]$-modules. Suppose that
we are given a short exact sequences of $\lie a$-modules $$
\begin{CD}
0 @>>> V @>\iota>> U @>\pi>> W @>>>0
\end{CD}
$$ and suppose that $V=V(\rho_V,\eta_V)$, $W=W(\rho_W,\eta_W)$ are
$\lie a[t]$-modules. A natural question is whether $U$ admits a
$\lie a[t]$-module structure with makes the above short exact
sequence into an extension of $\lie a[t]$-modules. A priori, for a
general extension of $\lie a$-modules, there is no reason to
expect that such a structure exists. However, if the corresponding
extension is equivalent to the split extension, then $U$ admits a
canonical $\lie a[t]$-module structure corresponding to the direct
sum of $V$ and $W$ as $\lie a[t]$-modules. In this section we
classify all possible $\lie a[t]$-module extensions as
deformations of that canonical one. This describes, in particular,
all extensions of $\lie a[t]$-modules such that the corresponding
category of $\lie a$-modules is semisimple.

We begin with some standard
results on $\Ext$.

\subl{Ex0} Recall that a triple~$(U,\iota,\pi)$, where~$U$ is a
$\lie a[t]$-module, $\iota\in\Hom_{\lie a[t]}(W,U)$, $\pi\in
\Hom_{\lie a[t]}(U,V)$ and
$$
\begin{CD}
0 @>>> W @>\iota>> U @>\pi>> V @>>> 0
\end{CD}
$$
is a short exact sequence of~$\lie a[t]$-modules, is called an
extension of~$W$ by~$V$. Two such triples~$(U,\iota,\pi)$
and~$(U',\iota', \pi')$ are equivalent if and only if there exists
an isomorphism of~$\lie a[t]$-modules $\psi:U\to U'$ such
that~$\psi\circ\iota=\iota'$ and~$\pi'\circ\psi=\pi$. It is
well-known that~$\Ext^1_{\lie a[t]}(V,W)$ identifies with the set
of equivalence classes of extensions of~$W$ by~$V$ and we
write~$[U,\iota,\pi]$ for the class of~$(U,\iota,\pi)$.

The set~$\Ext^1_{\lie a[t]}(V,W)$ can be endowed with a structure
of a vector space in the following way. Let~$[U_k,\iota_k,\pi_k]$,
$k=1,2$ be two elements of~$\Ext^1_{\lie a[t]}(V,W)$.
Then~$[U_1,\iota_1,\pi_1]+[U_2,\iota_2,\pi_2]=[U, \iota,\pi]$
where~$U$, $\iota$ and~$\pi$ are defined in the following way.
Let~$\wdt{U}= \{(w,u_1,u_2)\,:\, w\in W, u_i\in U_i,
\pi_1(u_1)=\pi_2(u_2)\}$ which is a submodule of~$W\oplus
U_1\oplus U_2$. Let $N$ be a submodule of~$\wdt{U}$ consisting of
all elements of the form~$(w_1+w_2,-\iota_1(w_1),-\iota_2(w_2))$.
Let~$U=\wdt{U}/N$ and define~$\iota(w)=(w,0,0)+N$, $\pi(
(w,u_1,u_2)+N)=\pi_1(u_1)= \pi_2(u_2)$. The extension
$(U,\iota,\pi)$  is called the Baer sum of
extensions~$(U_1,\iota_1,\pi_1)$ and~$(U_2,\iota_2,\pi_2)$.

The multiplication by~$z\in\bc$ in~$\Ext^1_{\lie a[t]}(V,W)$ is
defined as follows. Given~$[U,\iota,\pi]$, define~$z[U,\iota,\pi]$
to be the equivalence class of~$[U',\iota',\pi']$ where
$U'=(W\oplus U)/M$ with~$M=\{(-z w, \iota(w)): w\in W\}$,
$\iota'(w)=(w,0)+M$ and $\pi'((w,u)+M)=\pi(u)$.

Finally, the zero element of~$\Ext^1_{\lie a[t]}(V,W)$ is the
equivalence class of the split extension~$(V\oplus
W,\iota_0,\pi_0)$ where~$\iota_0(w)=(0,w)$ and~$\pi_0((v,w))=v$.

\subl{Ext10} Given $(\rho_i,\eta_i)\in\cal{H}(V_i)$,  $i=1,2$ let
$(\rho,\eta)=(\rho_1\oplus\rho_2,\eta_1\oplus\eta_2)$ be the
obvious maps from $\lie a\to\End(V_1\oplus V_2)$. Moreover, if
$\tilde\eta:\lie a\to\Hom(V_1,V_2)$ is any linear map, we regard
it as a map $\lie a\to \End(V_1\oplus V_2)$ also in the obvious
way. Let $$ \cal{E}(V_1(\rho_1,\eta_1),
V_2(\rho_2,\eta_2))=\{\tilde\eta\in \Hom(\lie a, \Hom(V_1,
V_2):(\rho,\eta+\tilde\eta)\in\cal{H}(V_1\oplus V_2)\}. $$ It is
obvious that
$$
\cal E(V_1(\rho_1,\eta_1),V_2(\rho_2,\eta_2))\subset \Hom_{\lie
a}(\lie a, \Hom(V_1,V_2)),
$$
where $\Hom(V_1,V_2)$ is regarded as an $\lie a$-module in a
natural way. We claim that $\cal
E(V_1(\rho_1,\eta_1),V_2(\rho_2,\eta_2))$ is a vector subspace of
$\Hom_{\lie a}(\lie a,\Hom(V_1,V_2))$. Indeed, since $\eta+\tilde
\eta\in\cal H_{\rho}(V_1\oplus V_2)$, we conclude that
$\ker(\eta+\tilde \eta)_T\supset \ker\tau[t]$. Now, since
$\tilde\eta(V_2)=0$ and $\tilde\eta(V_1)=0$, we have
$\tilde\eta(x)\circ \psi\circ\tilde\eta(y)=0$ for all $x,y\in\lie
a$ and for all $\psi\in\End(V_2)$. It follows immediately that
$$
(\eta+\tilde \eta)_T|_{T^r(\lie a)}=\eta^{\odot
r}+\sum_{j=0}^{r-1} \eta_2^{\odot j}\odot\tilde\eta \odot
\eta_1^{\odot r-j-1}.
$$
Since~$\eta\in\cal H_{\rho}(V)$, $\eta^{\odot r}(\ker\tau_r)=0$.
Thus, $\eta+\tilde\eta\in\cal H_{\rho}(V)$ if and only if
\begin{enumerit}
\item[$(E_1)$] $\tilde\eta\in\Hom_{\lie a}(\lie a,\Hom(V_1,V_2))$,
\item[$(E_2)$] $\sum_{j=0}^{r-1} \eta_2^{\odot j} \odot
\tilde\eta\odot \eta_1^{\odot r-j-1}\big( \ker\tau_r\big)=0$, $r>1$.
\end{enumerit}
In particular, $\mathcal E(V_1(\rho_1,\eta_1),V_2(\rho_2,\eta_2))$
is a vector subspace of $\Hom_{\lie a}(\lie a,\Hom(V_1,V_2))$.

\subl{Ext10a} 
Given $\psi\in\Hom_{\lie a}(V_1,V_2)$, set $\tilde\eta_\psi(x):=\eta_2(x)\circ\psi-\psi\circ\eta_1(x)\in
\Hom(\lie a,\Hom(V_1,V_2))$ and define
$$
\mathcal E_0(V_1(\rho_1,\eta_1),V_2(\rho_2,\eta_2))=\{\tilde\eta_\psi\,:\,\psi\in\Hom_{\lie a}(V_1,V_2)\}.
$$
\begin{lem} The set $ \mathcal E_0(V_1(\rho_1,\eta_1),V_2(\rho_2,\eta_2))$
is a subspace of $\Hom_{\lie a}(\lie a,\Hom(V_1,V_2))$. Further,
the assignment  $\psi\mapsto\tilde\eta_\psi$ defines an
isomorphism of vector spaces, $$ 
\mathcal
E_0(V_1(\rho_1,\eta_1),V_2(\rho_2,\eta_2))\cong \Hom_{\lie
a}(V_1,V_2)/ \Hom_{\lie a[t]}(V_1,V_2).$$
\end{lem}
\begin{pf} The first statement of the lemma is established in the
course of this section where it is shown that  $ \mathcal
E_0(V_1(\rho_1,\eta_1),V_2(\rho_2,\eta_2))$ is the kernel of a
linear map. The linear map $\psi\mapsto\tilde\eta_\psi$ is surjective
by definition and its kernel is $\Hom_{\lie a[t]}(V_1,V_2)$ by~\lemref{Cons22}.
\end{pf}

\subl{Ext10b} The following theorem is the main result of this section.
\begin{thm}\label{extthm}
Suppose that $\Ext^1_{\lie a}(V_1,V_2)=0$. Then we have an
isomorphism of vector spaces $$ 
\be:\mathcal
E(V_1(\rho_1,\eta_1),V_2(\rho_2,\eta_2))/ \mathcal
E_0(V_1(\rho_1,\eta_1),V_2(\rho_2,\eta_2))\stackrel{\sim}{\longrightarrow} \Ext^1_{\lie
a[t]}(V_1(\rho_1,\eta_1),V_2(\rho_2,\eta_2)),$$ 
given by the assignment~$\be:[\tilde\eta]\mapsto [
V(\rho,\eta+\tilde\eta),\iota_0,\pi_0]$, where~$[\tilde\eta]$
denotes the class of~$\tilde\eta$ modulo~$\mathcal E_0$,
$V=V_1\oplus V_2$, $\rho=\rho_1\oplus\rho_2$,
$\eta=\eta_1\oplus\eta_2$, $\iota_0(v_2)=(0,v_2)$
and~$\pi_0((v_1,v_2))=v_1$ for all~$v_i\in V_i$.
\end{thm}
\noindent
The proof of this theorem is given in~\ref{Ext12}--\ref{Ext17}.

\subl{Ext12}
Throughout~\ref{Ext12}-\ref{Ext17} we keep the notations of~\ref{Ext10}.
The first step of our proof is
\begin{lem}
We have the following short exact sequence of
modules:
$$
\begin{CD}
0 @>>> V_2(\rho_2,\eta_2) @>\iota_0>> V
(\rho,\eta+\tilde\eta) @>\pi_0>> V_1(\rho_1,\eta_1)@>>> 0.
\end{CD}
$$
\end{lem}
\begin{pf}
We only need to check that $\iota_0$ and~$\pi_0$ are homomorphisms
of $\lie a[t]$-modules, since the above sequence is obviously
exact as a sequence of $\lie a$-modules. 
Indeed, we have for all~$x\in\lie a$,
\begin{align*}
&\iota_0( \eta_2(x)(v_2))=(0,\eta_2(x)(v_2))=(\eta+\tilde\eta)(x)((0,v_2))=
(\eta+\tilde\eta)(x)(\iota_0(v_2)),\\
&\pi_0 ( (\eta+\tilde\eta)(x)((v_1,v_2)))=(\eta_1(x)(v_1),\eta_2(x)(v_2)+
\tilde\eta(x)(v_1))=\eta_1(x)(v_1)=\eta_1(\pi_0((v_1,v_2))).\qedhere
\end{align*}
\end{pf}

\subl{Ext15}
Lemma~\ref{Ext12} allows us to define
a map
$$
\begin{array}{rcl}
\tilde{\be}:\cal{E}(V_1(\rho_1,\eta_1), V_2(\rho_2,\eta_2))&\longrightarrow&
\Ext^1_{\lie a[t]}(V_1(\rho_1,\eta_1), V_2(\rho_2,\eta_2))\\
{}[\tilde{\eta}]&\longmapsto&[ V(\rho,\eta+\tilde\eta)].
\end{array}
$$
\begin{lem}
\begin{enumerit}
\item $\tilde{\be}$ is a homomorphism of vector spaces.
\item $\ker\tilde{\be}=\mathcal{E}_0(V_1(\rho_1,\eta_1),V_2(\rho_2,\eta_2))$.
In particular, $\cal E_0(V_1(\rho_1,\eta_1),V_2(\rho_2,\eta_2))$ is a
vector subspace of $\cal E(V_1(\rho_1,\eta_1),V_2(\rho_2,\eta_2))$.
\end{enumerit}
\end{lem}
\begin{pf}
The first step is to check that the Baer sum of~$[V(\rho,\eta+\tilde\eta_1),
\iota_0,\pi_0]$ and~$[V(\rho,\eta+\tilde\eta_2), \iota_0,\pi_0]$ equals
$[V(\rho,\eta+\tilde\eta_1+\tilde\eta_2),\iota_0,\pi_0]$. Indeed, by~\ref{Ex0}
$[V(\rho,\eta+\tilde\eta),\iota_0,\pi_0]+[V(\rho,\eta+\tilde\eta),
\iota_0,\pi_0]$ is defined as~$[U,\iota,\pi]$ where $U=\wdt{U}/N$
with
$$
\wdt{U}=\{(w, (v,w_1), (v,w_2))\,:\, v\in V_1,\, w,w_1,w_2\in
V_2\}
$$
and
$$
N=\{ (w_1+w_2,(0,-w_1),(0,-w_2))\,:\, w_1,w_2\in
V_2\}\subset\wdt{U}.
$$
The embedding~$\iota:V_2(\rho_2,\eta_2)\to U$ and the
projection~$\pi:U\to V_1(\rho_1,\eta_1)$ are defined
by
$$
\iota(w)=(w,(0,0),(0,0))+N,\qquad \pi( (w,(v,w_1),(v,w_2))+N)=v.
$$
Observe also that $U=U(\rho_U,\eta_U)$ where
$$
\rho_U(x)( (w,(v,w_1),(v,w_2))+N)=
(\rho_2(x)(w),\rho(x)((v,w_1)),
\rho(x)((v,w_2)))+N
$$
and
$$
\eta_U(x)((w,(v,w_1),(v,w_2))+N)=
(\eta_2(x)(w),(\eta+\tilde\eta_1)(x)((v,w_1)),
(\eta+\tilde\eta_2)(x)((v,w_2)))+N.
$$
Define~$\Psi:U\to V(\rho,\eta+\tilde\eta_1+\tilde\eta_2)$ by
$\Psi( (w,(v,w_1),(v,w_2))+N)=(v,w+w_1+w_2)$. This map is
obviously well-defined, satisfies $\Psi\circ\iota=\iota_0$, $\pi_0\circ\Psi=\pi$ and is an isomorphism of $\lie a$-modules.
It remains to check that~$\Psi$ is an isomorphism of~$\lie a[t]$-modules. We have
\begin{align*}
\Psi\circ\eta_U(x)((w,(v,w_1),(v,w_2))+N)&=
(\eta_1(x)(v),\eta_2(x)(w+w_1+w_2)+(\tilde\eta_1+\tilde\eta_2)(x)(v))\\
&=(\eta+\tilde\eta_1+\tilde\eta_2)(x)( (v,w+w_1+w_2))\\& =
(\eta+\tilde\eta_1+\tilde\eta_2)(x)\circ \Psi( (w,(v,w_1),(v,w_2))+N).
\end{align*}
Thus, 
$\Psi\circ\eta_U(x)=
(\eta+\tilde\eta_1+\tilde\eta_2)(x)\circ\Psi$, and so~$\Psi:U\to
V(\rho, \eta+\tilde\eta_1+\tilde\eta_2)$ is an isomorphism
of~$\lie a[t]$-modules by~\lemref{Cons22}.

Furthermore, let~$z\in\bc$. We claim that $z[V(\rho,\eta+\tilde\eta)
\iota_0,\pi_0]=[V(\rho,\eta+z\tilde\eta),\iota_0,\pi_0]$. Indeed,
by~\ref{Ex0}, $z[V(\rho,\eta+\tilde\eta),\iota_0,\pi_0]=[U',\iota',\pi']$
where
$$
\begin{array}{lll}
U'=W\oplus V(\rho,\eta+\tilde\eta)/M,\qquad\qquad &M=\{ (-z w,(0,w))\,:\, w\in
V_2\},
\\
\iota'(w)=(w,(0,0))+M, & \pi'( (w,(v,w'))+M )=v,
\qquad &\forall\, v\in V_1,\, w,w'\in V_2
\end{array}
$$
Define~$\Psi':U'\to V$ as~$\Psi'( (w,(v,w'))+M)=(v,w+zw')$.
Then~$\Psi'$ is a well-defined isomorphism of $\lie a$-modules and
it is easy to check that~$\Psi'\circ \iota'=\iota_0$,
$\pi_0\circ\Psi'=\pi'$. Observe that~$U'=U'(\rho',\eta')$, where
$$
\rho'(x)( (w,(v,w'))+M)=(\rho_2(x)(w),\rho(x)(v,w'))+M
$$
and
$$
\eta'(x)( (w,(v,w'))+M)=(\eta_2(x)(w),(\eta+\tilde\eta)(v,w'))+M.
$$
We have
\begin{align*}
\Psi'\circ\eta'(x)( (w,(v,w'))+M)&=
\Psi'((\eta_2(x)(w), (\eta_1(x)(v),\eta_2(x)(w')+\tilde\eta(x)(v)))+M)\\
&=( \eta_1(x)(v),\eta_2(x)(w+zw')+z\tilde\eta(x)(v))=
(\eta+z\tilde\eta)(x)((v,w+zw'))\\&=
(\eta+z\tilde\eta)(x)\circ \Psi'( (w,(v,w'))+M).
\end{align*}
Therefore~$\Psi':U\to V(\rho,\eta+z\tilde\eta)$ is an isomorphism
of~$\lie a[t]$-modules by~\lemref{Cons22} and
so~$[U',\iota',\pi']=[V(\rho,\eta+z\tilde\eta), \iota_0,\pi_0]$.

It remains to prove~(ii). Suppose~$\tilde\eta\in\ker\wdt{\be}$.
Then $[V(\rho,\eta+\tilde\eta), \iota_0,\pi_0]=[V(\rho,\eta),\iota_0,\pi_0]$.
Let~$\Psi_0: V(\rho,\eta+\tilde\eta)\to
V(\rho,\eta)$ be a corresponding isomorphism
of~$\lie a[t]$-modules. We may write $ \Psi_0(
(v_1,v_2))=(\psi_{1,1}(v_1)+\psi_{1,2}(v_2),\psi_{2,1}(v_1)+\psi_{2,2}(v_2))$.
Since~$\Psi_0\circ\iota_0=\iota_0$ and~$\pi_0\circ\Psi_0=\pi_0$ we
conclude that~$\psi_{1,2}=0$, $\psi_{1,1}=\id_{V_1}$,
$\psi_{2,2}=\id_{V_2}$ and~$\psi_{2,1}=\psi\in\Hom(V_1,V_2)$.
Since $\Psi_0$ is an isomorphism of $\lie a$-modules, it follows
that $\psi\in\Hom_{\lie a}(V_1,V_2)$. Furthermore,
$\Psi_0\circ(\eta+\tilde\eta)(x)=\eta(x)\circ\Psi_0$
by~\lemref{Cons22}. Then
\begin{align*}
\Psi_0\circ (\eta+\tilde\eta)(x)( (v_1,v_2))&= \Psi_0(
\eta_1(x)(v_1),\eta_2(x)(v_2)+\tilde\eta(x)(v_1))\\&=
(\eta_1(x)(v_1),\eta_2(x)(v_2)+\tilde\eta(x)(v_1)+\psi\circ\eta_1(x)(v_1))\\
&=\eta(x)\circ\Psi_0((v_1,v_2))=(\eta_1(x)(v_1),\eta_2(x)(v_2)+
\eta_2(x)\circ\psi(v_1)).
\end{align*}
It follows
that~$\tilde\eta=\tilde\eta_\psi$ and so $\ker\tilde{\be}\subset\mathcal E_0(V_1(\rho_1,\eta_1),
V_2(\rho_2,\eta_2))$. For the opposite inclusion,
let~$\tilde\eta\in \mathcal
E_0(V_1(\rho_1,\eta_1),V_2(\rho_2,\eta_2))$. Then $\tilde\eta=\tilde\eta_\psi$
for some~$\psi\in\Hom_{\lie a}(V_1,V_2)$. 
Define~$\Psi_0:
V(\rho,\eta+\tilde\eta) \to V_1(\rho_1,\eta_1)\oplus
V_2(\rho_2,\eta_2)$ by~$\Psi_0((v_1,v_2))= (v_1,v_2+\psi(v_1))$.
One checks as above that~$\Psi_0$ is an isomorphism of $\lie
a[t]$-modules, $\Psi_0\circ\iota_0=\iota_0$,
$\pi_0\circ\Psi_0=\pi_0$. Thus, $\wdt{\be}(\tilde\eta)$ is
the class of the split extension and so~$\tilde\eta\in\ker\wdt{\be}$.
\end{pf}

\subl{Ext17}
The following lemma completes the proof of our theorem.
\begin{lem}
The homomorphism~$\tilde{\be}$ defined in~\ref{Ext15} is surjective.
\end{lem}
\begin{pf}
Let $(U,\iota,\pi)$ be an extension of $V_2(\rho_2,\eta_2)$ by
$V_1(\rho_1,\eta_1)$. Observe that~$U=U(\rho_U,\eta_U)$ for some
pair $(\rho_U,\eta_U)\in\mathcal H(U)$. Since $\Ext^1_{\lie
a}(V_1,V_2)=0$, $(U,\iota,\pi)$ splits as an extension of $\lie
a$-modules. Then there exists $j\in\Hom(V_1,U)$ such that
$U=j(V_1)\oplus \iota(V_2)$, $\pi\circ j=\id$ and $\rho_U(x)\circ
j=j\circ\rho_1(x)$. Observe that, since~$\pi$, $\iota$ are
homomorphisms of $\lie a[t]$-modules,
$\pi\circ\eta_U(x)=\eta_1(x)\circ\pi$, $\eta_U(x)\circ\iota=
\iota\circ\eta_2(x)$ for all~$x\in\lie a$.
One has
$$
\pi\circ(\eta_U(x)\circ j-j\circ \eta_1(x))=
\eta_1(x)-\eta_1(x)=0.
$$
Therefore, $\eta_U(x)\circ j-j\circ\eta_1(x)\in \Hom(\lie
a,\Hom(V_1,\ker\pi))=\Hom(\lie a,\Hom(V_1,\Im\iota))$. Thus,
$$
\tilde\eta_j(x):=\iota^{-1}(\eta_U(x)\circ j-j\circ\eta_1(x))
$$
is a well-defined element of~$\Hom(\lie a,\Hom(V_1,V_2))$.

Let us check that~$\tilde\eta_j\in\mathcal E(V_1(\rho_1,\eta_1),
V_2(\rho_2,\eta_2))$. 
By construction, $\tilde\eta_j$ is a
homomorphism of $\lie a$-modules. Furthermore,
\begin{align*}
\sum_{s=0}^{r-1} \eta_2^{\odot s}\odot \tilde\eta_j\odot \eta_1^{\odot r-s-1}&=
\iota^{-1}\Big(\sum_{s=0}^{r-1} \eta_U^{\odot s+1}\circ j\circ
\eta_1^{\odot r-s-1}-\sum_{s=0}^{r-1} \eta_U^{\odot s}\circ j\circ
\eta_1^{\odot r-s}\Big)\\
&=\iota^{-1}( \eta_U^{\odot r}\circ j-j\circ \eta_1^{\odot r}).
\end{align*}
Since both~$\eta_U$ and~$\eta_1$ satisfy~$(C_3)$,
it follows that~$\tilde\eta_j$ satisfies~$(E_2)$.

It remains to prove that~$[U,\iota,\pi]=[V(\rho,\eta+
\tilde\eta_j),\iota_0,\pi_0]=\tilde{\be}(\tilde\eta_j)$. For,
define~$\Psi: V(\rho,\eta+\tilde\eta_j)\to U$ by~$\Psi(
(v_1,v_2))=j(v_1)+\iota(v_2)$. Evidently, $\Psi$ is an isomorphism
of $\lie a$-modules, $\Psi\circ\iota_0=\iota$ and
$\pi\circ\Psi=\pi_0$. One has 
\begin{multline*}
\Psi\circ(\eta+\tilde\eta_j(x))((v_1,v_2))=
j(\eta_1(x)(v_1))+\iota(\eta_2(x)(v_2))+\eta_U(x)(j(v_1))-
j(\eta_1(x)(v_1))\\=\eta_U(x)(j(v_1)+\iota(v_2))=
\eta_U(x)\circ\Psi((v_1,v_2)).
\end{multline*}
It follows from~\lemref{Cons22} that~$\Psi$ is an isomorphism of
$\lie a[t]$-modules.
\end{pf}

\subl{Ext40} Assume from now on that $\lie g$ is a simple Lie algebra,
$R^+$ a set of positive roots for $\lie g$ with respect to a
Cartan subalgebra $\lie h$ of $\lie g$.  Let $\lie g=\lie
n^+\oplus\lie h\oplus\lie n^-$ be the corresponding triangular
decomposition. Let $\theta$ be the highest root of $\lie g$.

Let $\cal F$ be the category of finite-dimensional $\lie
g[t]$-modules. We shall prove the following in the rest of the section.
\begin{thm} \label{extsimple}
Let $V\in\cal F$ be irreducible. Then,
$$
\dim\Ext_{\lie g[t]}^1(\bc, V)\le 1,
$$ with  equality holding if and only if $V\cong
\lie g(\ad, a\ad)$ for some $a\in\bc$. Moreover, if $V'\in\cal F$ is
also irreducible, then
$$\dim\Ext^1_{\lie g[t]}(V,V')=\sum_{a\in\bc} \dim \Hom_{\lie g[t]}(\lie g(\ad,
a\ad), V^*\otimes V').
$$
In particular, $\Ext^1_{\lie g[t]}(V,V)\not=0$.
\end{thm}

\subl{Ext50} We shall need the following proposition on the structure
of $\cal F$.
Recall that a module in $\cal F$ is completely
reducible if and only if it is isomorphic to a direct sum of irreducible
$\lie g[t]$-modules.
\begin{prop} 
\begin{enumerit}
\item\itlabel{Ext50.i}  Let $V\in \cal F$ be irreducible. Then either $V$ is
trivial or there exist
$k\in\bn$, simple non-trivial finite-dimensional $\lie g$-modules $V_i(\rho_i)$ and distinct
$a_i\in\bc$, $1\le i\le k$ such that
$$
V\cong V_1(\rho_1,a_1\rho_1)\otimes\cdots\otimes V_k(\rho_k,a_k\rho_k).
$$
In particular there exists a unique element $v\in V$ such that
$V=U(\lie g[t])v$, $\lie n^+[t]v=0$ and $\lie h[t]v=\lie
hv\subset\bc v$. Moreover if $V$ is non-trivial, then $\lie
hv=\bc v$.

\item\itlabel{Ext50.ii}
Given $k\in\bn$ and   irreducible $\lie g$-modules
$V_i(\rho_i)$, $1\le i\le k$, the $\lie g[t]$-module
$V_1(\rho_1,a_1\rho_1)\otimes\cdots\otimes V_k(\rho_k,a_k\rho_k)$ is
irreducible if and only if  $a_i\ne a_j$, $1\le i\le j\le k$ and is
completely reducible otherwise.

\item\itlabel{Ext50.iii}
Let $V\in\cal F$. Then $V$ is completely reducible if and only if  there exists $k\in\bn$ and distinct
$a_i\in\bc$, $1\le i\le k$ such that
$$
(x\otimes (t-a_1)\cdots (t-a_k))V=0,
$$
for all $x\in\lie g$.
\item\itlabel{Ext50.iv} Let $V\in\cal F$ be non-trivial. There exists $a\in\bc$
such that $\Hom_{\lie g[t]}(\lie g(\ad, a\ad), V\otimes V^*)\ne 0$.

\end{enumerit}
\end{prop}\label{stf}
\begin{pf} Parts \eqref{Ext50.i} and \eqref{Ext50.ii} are well-known (cf.~\cite{C,CPnew} for instance).
For \eqref{Ext50.iii}, notice that if
$V\cong V_1(\rho_1,a_1\rho_1)\otimes\cdots\otimes
V_k(\rho_k,a_k\rho_k)$ is  then
$$
(x\otimes (t-a_1)\cdots (t-a_k))V=0,
$$
for all $x\in\lie g$ and hence the result follows
for all completely reducible modules in $\cal F$.
For the
converse statement, consider the map $\ev_{a_1,\dots ,a_k}:\lie
g[t]\to\lie g^{\oplus k}$ defined by
$$
\ev_{a_1,\dots,a_k}(x\otimes t^s)=(a_1^sx,\dots,a_k^sx), \qquad x\in\lie g.
$$
This map is a homomorphism of Lie algebras and is surjective if
the $a_i$ are distinct for $1\le i\le k$.
Thus the action of $\lie g[t]$ on  $V$ gives rise to an action of  $\lie g^{\oplus k}$.
Since $\lie g^{\oplus k}$ is a semisimple Lie algebra, $V$ is
completely reducible as a $\lie g^{\oplus k}$-module. Each
irreducible component is in fact  an irreducible $\lie
g[t]$-module obtained by pulling back by~$\ev_{a_1,\dots,a_k}$ and hence~\eqref{Ext50.iii} follows.

Clearly, it is sufficient to prove~\eqref{Ext50.iv} for $V$ simple and non-trivial.
Observe that $V(\rho,a\rho)^*\cong V^*(\rho^*,a\rho^*)$ as $\lie g[t]$-modules, whence 
$$\Hom_{\lie g[t]}(\bc, V(\rho,a\rho)\otimes
V^*(\rho^*,a\rho^*) )\ne 0$$ 
and  also that 
$$\Hom_{\lie g[t]}(\lie g(\ad,a\ad), V(\rho,a\rho) \otimes V^*(\rho^*,a\rho^*))\ne
0.
$$
It follows immediately that if $V\cong V_1(\rho_1,a_1\rho_1)\otimes\cdots\otimes V_k(\rho_k,a\rho_k)$,
then $\Hom_{\lie g[t]}(\lie g(\ad,a\ad), V\otimes V^*)\ne 0$.
\end{pf}

\noindent
Let $x^-_\theta\in\lie g$ be a non-zero element in the root space
corresponding to the root $-\theta$.
\begin{cor} 
Let $V\in\cal F$. Then $V$
is completely reducible if and only if  there exists $k,\ell\in\bn$, 
distinct $a_i\in\bc$, $1\le i\le k$ and elements $v_j\in V$, $1\le
j\le \ell$ such that $\lie n^+[t]v_j=0$, $\lie h[t]v_j\subset\bc
v_j$ and
$$
(x^-_\theta\otimes (t-a_1)\cdots (t-a_k))v_j=0,  \qquad 1\le j\le \ell.
$$
\end{cor}
\begin{pf} The corollary is immediate from the proposition if we
prove that the condition
$$
(x^-_\theta\otimes (t-a_1)\cdots
(t-a_k))v=0,
$$
for some $v\in V$ with $\lie n^+[t]v=0$, $\lie
h[t]v\subset\bc v$ implies  that
$$
(x\otimes (t-a_1)\cdots (t-a_k)) U(\lie g[t])v=0,
$$
for all $x \in\lie g$. Since
$U(\lie g[t])v=U(\lie n^-[t])v$, and $[x^-_\theta,\lie
n^-]=0$, we see that
$$(x^-_\theta\otimes (t-a_1)\cdots (t-a_k))U(\lie
g[t])v=0.
$$
The result follows now since any $x\in\lie g$ can be written as a linear combination
of elements of the form $[y_1,[y_2,\dots[y_m,x^-_\theta]]\cdots]$
for $m\in\bn$ and $y_1,\dots ,y_m\in\lie g$.
\end{pf}

\subl{Ext60}
\begin{proof}[Proof of~\thmref{extsimple}]
Suppose first
that  that~$V\cong \lie g(\ad,a\ad)$ for some~$a\in\bc$. Then $\Hom_{\lie g}(\bc,V)=0$
and so~$\cal E_0(\bc,V)=0$. Thus,
$\dim\Ext^1_{\lie g[t]}(\bc,V)\le\dim\Hom_{\lie g}(\lie g,\Hom(\bc,V))=1$.
Since the $\lie g[t]$-module described in~\ref{Cons35} provides a
non-zero element of $\Ext_{\lie g[t]}^1(\bc,V)$, it follows that
$\dim\Ext_{\lie g[t]}^1(\bc,V)=1$.

Now suppose that $V\cong V_1(\rho_1,a_1\rho_1)\otimes\cdots\otimes
V_k(\rho_k,a_k\rho_k)$ is not isomorphic to $\lie g(\ad,a\ad)$. If $V$
is trivial, then the result is obvious. Otherwise, assume that we
have a short exact sequence of $\lie g[t]$-modules,
\begin{equation}\label{Ext50.10}
\begin{CD}
0@>>> \bc @>\iota>> M @>\pi>> V @>>> 0.
\end{CD}
\end{equation}
Let $m\in M$ be such that $\pi(m)\in V$
is the element $v$ of \propref{Ext50}\eqref{Ext50.i}. Assume that
$hv=\lambda(h)v$ for some $\lambda\in\lie h^*$ and $hm=\lambda(h)m$
as well. Now, since
$$
0=x^-_\theta\otimes (t-a_1)\cdots(t-a_k)v=\pi\Big(x^-_\theta\otimes(t-a_1)\cdots
(t-a_k)m\Big),
$$
it follows that $(x^-_\theta\otimes (t-a_1)\cdots
(t-a_k))m\ne 0$ only if it lies in~$\Im\iota\cong\bc$, which implies $\lambda=\theta$.

Suppose first that~$\lambda\not=\theta$. Then $(x^-_\theta\otimes (t-a_1)\cdots
(t-a_k))U(\lie g[t])m= 0$. Since $M=V\oplus\bc$ as $\lie
g$-modules and $(x^-_\theta\otimes (t-a_1)\cdots
(t-a_k))U(\lie g[t])\bc= 0$ it follows that
$$
(x^-_\theta\otimes
(t-a_1)\cdots (t-a_k))M= 0,
$$
which implies by~\corref{Ext50}
that $M$ is completely reducible and hence the short
exact sequence~\eqref{Ext50.10} splits.

Suppose now that $\lambda=\theta$ and that
$V=V_1(\rho_1,a_1\rho_1)\otimes\cdots\otimes V_k(\rho_k,a_k\rho_k).$
Then it is easy to see from the representation theory of simple
Lie algebras that $k\le 2$. Since $k=1$ is the case when $V\cong
\lie g(\ad,a\ad)$, it is enough to consider the case $k=2$. But this can occur only if
$\lie g\cong\lie{sl}_{n+1}$ and only if one of the modules $V_i(\rho_i)$ is isomorphic
to the natural representation while the other is isomorphic to its
dual. Then
$$
V\cong\lie g\oplus\bc
$$
as $\lie g$-modules, whence
$$
\dim\Ext^1_{\lie g[t]}(\bc, V)\le\dim\Hom_{\lie g}(\lie g,\Hom(\bc,V))=
\dim\Hom_{\lie g}(\lie g,V)=1
$$
by~\thmref{extthm}.
Since $\dim\Hom_{\lie g}(\bc, V)=1$
and $\dim\Hom_{\lie g[t]}(\bc,V)=0$, it follows from~\lemref{Ext10a} that $\dim\cal E_0(\bc,V)=1$
and so $\dim\Ext^1_{\lie g[t]}(\bc,V)=0$.

The last statement is now an immediate consequence of
the fact that for all $V,V'\in\cal F$, we have
\begin{equation*}\Ext^1_{\lie g[t]}(V,V')\cong\Ext^1_{\lie g[t]}(\bc, V^*\otimes V').\qedhere
\end{equation*}
\end{proof}

\end{document}